\documentclass{amsart}

\usepackage{amssymb,amsfonts, latexsym,amsmath,hhline,array,longtable}
\usepackage{pdfsync,color, comment,colortbl, mathrsfs,stmaryrd,cite,graphicx}

\usepackage{ifpdf}
\ifpdf \usepackage[colorlinks=true, citecolor=blue, linkcolor=blue, urlcolor=blue]{hyperref} \fi

\newcommand{\cal}{\mathcal}

\newtheorem{formula}{}[section]
\newtheorem{definition}[formula]{Definition}
\newtheorem{corollary}[formula]{Corollary}
\newtheorem{remark}[formula]{Remark}
\newtheorem{lemma}[formula]{Lemma}
\newtheorem{theorem}[formula]{Theorem}

\def\thrm{\begin{theorem}}
\def\thrml#1{\begin{theorem}\label{#1}}
\def\ethrm{\end{theorem}}
\def\rmrk{\begin{remark}}
\def\rmrkl#1{\begin{remark}\label{#1}}
\def\ermrk{\end{remark}}
\def\dfntn{\begin{definition}}
\def\dfntnl#1{\begin{definition}\label{#1}}
\def\edfntn{\end{definition}}
\def\nmrt{\begin{enumerate}}
\def\enmrt{\end{enumerate}}
\def\tm#1{\item[{\rm (#1)}]}
\def\qtnl#1{\begin{equation}\label{#1}}
\def\eqtn{\end{equation}}
\def\lmm{\begin{lemma}}
\def\lmml#1{\begin{lemma}\label{#1}}
\def\elmm{\end{lemma}}
\def\crllr{\begin{corollary}}
\def\crllrl#1{\begin{corollary}\label{#1}}
\def\ecrllr{\end{corollary}}
\def\css{\begin{cases}}
\def\ecss{\end{cases}}
\def\prf{\begin{proof}}
\def\eprf{\end{proof}}

\def\cP{{\cal P}}
\def\cQ{{\cal Q}}

\def\cX{{\cal X}}

\def\mN{{\mathbb N}}

\def\fK{{\mathfrak K}}

\DeclareMathOperator{\aut}{Aut}

\DeclareMathOperator{\GaL}{{\rm \Gamma}L}

\DeclareMathOperator{\poly}{poly}

\DeclareMathOperator{\soc}{Soc}

\DeclareMathOperator{\sym}{Sym}

\def\qaq{\quad\text{and}\quad}

\begin{document}

\title{On computing the closures of solvable permutation groups}
\author{Ilia\ Ponomarenko}
\address{Steklov Institute of Mathematics at St.\ Petersburg,  Russia}
\email{inp@pdmi.ras.ru}
\author{Andrey\ V.\ Vasil'ev}
\address{Sobolev Institute of Mathematics, Novosibirsk, Russia;}
\address{Novosibirsk State University, Novosibirsk, Russia}
\email{vasand@math.nsc.ru}
\thanks{The second author was supported by the Program of Fundamental Research RAS, project FWNF-2022-0002.}
\date{}

\begin{abstract}
Let $m\ge 3$ be an integer. It is proved that the $m$-closure of a given solvable permutation group
of degree $n$ can be constructed in time $n^{O(m)}$.

\end{abstract}

\maketitle

\section{Introduction}

Let $m$ be a positive integer and let $\Omega$ be a finite set. The {\it $m$-closure} $G^{(m)}$ of $G\le\sym(\Omega)$ is the largest permutation group on $\Omega$ having the same orbits as $G$ in
its induced action on the cartesian power~$\Omega^m$. The $m$-closure of a permutation group was introduced by H.~Wielandt in~\cite{Wie1969}, where it was, in particular,  proved that $G^{(m)}$ can be treated as the full automorphism group of the set of all $m$-ary relations invariant with respect to~$G$. Since then the theory was developed in different directions, e.g., there were studied the closures of primitive groups~\cite{Liebeck1988a,PS1992,Yu2019}, the behavior of the closure with respect to permutation group operations~\cite{KaluK1976,Churikov2021,Vasilev2021}, totally closed abstract groups~\cite{AbdAT2022,AIPT2021,ChurPr2021}, etc.

From the computational point of view, the $m$-closure problem consisting in finding the $m$-closure of a given permutation group is of special interest; here and below, it is assumed that permutation groups are given by  generating sets, see \cite{Sere2002}. When the number~$m$ is given as a part of input, the problem seems to be very hard even if the input group is abelian. It is quite natural therefore to restrict the $m$-closure problem to the case when $m$ is fixed and the input group belongs to a certain class of groups. In this setting, polynomial-time algorithms for finding the $m$-closure were constructed for  the groups of odd order \cite{EvdP2001b} and, if $m=2$, also for nilpotent and supersolvable groups~\cite{Ponomarenko1994,Ponomarenko2020b}. Note that the case $m=1$ is trivial, because the $1$-closure of any permutation group~$G$ is equal to the direct product of symmetric groups acting on the orbits of~$G$.

The goal of the present paper is to solve the $m$-closure problem for $m\ge 3$ in the class of all solvable groups (note that  there is an efficient algorithm testing whether or not a given permutation group is solvable).

\thrml{040123b}
Given an integer $m\ge 3$,  the $m$-closure of a solvable permutation group of degree~$n$ can be found in time $n^{O(m)}$.
\ethrm

The proof of Theorem~\ref{040123b} is given in Section~\ref{040123b0}. A starting point in our approach to the proof is the main result in~\cite{OBrien2020} stating that for $m\ge 3$ the $m$-closure of every solvable permutation group is solvable. To apply this result, it suffices for a given solvable group~$G$ to find a solvable overgroup and then find $G^{(m)}$ inside it with the help of the Babai-Luks algorithm~\cite{BabaL1983}; the latter enables, in particular, to find efficiently the relative $m$-closure $G^{(m)}\cap H$ of an arbitrary group $G$ with respect to a solvable group $H$.

To explain  how to find the overgroup, we recall that a permutation group is said to be \emph{non-basic} if it is contained in a wreath product with the product action; it is \emph{basic} otherwise, see~\cite[Section~4.3]{Cam1999}. A classification of the primitive solvable linear groups having a faithful regular orbit~\cite{Yang2020} implies that for a sufficiently large primitive basic solvable group $G$, we have $G=G^{(m)}$ for all $m\ge 3$. This reduces the problem to solvable groups that are not basic, that is, to those that can be embedded in a direct or wreath product of smaller groups. In fact,  we only need to test whether the corresponding embedding exists and (if so) to find it explicitly. This is a subject of Section~\ref{010223a}.

All undefined terms can be found in \cite{Cam1999} (for permutation groups) and~\cite{Sere2002} (for permutation group algorithms).

The authors thank S.~V.~Skresanov  for fruitful discussions and useful comments.\medskip

\section{The embedding problem}\label{010223a}

Given permutation  groups $K\le\sym(\Delta)$ and $L\le\sym(\Gamma)$, we denote by $K\times L$ (respectively, $K\wr L$, $K\uparrow L$) the permutation group induced by the action of direct (respectively, wreath) product of~$K$ and~$L$ on  $\Delta\cup\Gamma$ (respectively, $\Delta\times\Gamma$, $\Delta^{|\Gamma|}$).

\thrml{090123a}
Let $m\ge 2$ be an integer, $K,L$ permutation groups, and $\star\in\{\times,\wr,\uparrow\}$. Then
$$
(K\star L)^{(m)}\le K^{(m)}\star L^{(m)}
$$
unless $\star=\uparrow$, $m=2$, and $K$ is $2$-transitive.
\ethrm
\prf
See \cite[Theorems~3.1, 3.2]{OBrien2020} and \cite[Theorem~1.2]{Vasilev2021}.
\eprf

Theorem~\ref{090123a}  is used to reduce the study of the $m$-closure of  a group $G\le\sym(\Omega)$ to permutation groups on smaller sets. From the algorithmic point of view, we need to solve the \emph{$\star$-embedding problem}: test whether there exists an embedding of~$G$ to~$K\star L$ for some sections~$K\le\sym(\Delta)$ and~$L\le\sym(\Gamma)$ of the group~$G$, such that $|\Delta|<|\Omega|$ and $|\Gamma|<|\Omega|$, and if so, then to find the embedding explicitly. By this, we mean finding a bijection~$f$ from~$\Omega$ to the underlying set of $K\star L$, such that
\qtnl{010223e}
f^{-1}Gf\le K\star L.
\eqtn
The $\star$-embedding problem is easy if $G$ is intransitive and $\star=\times$, or imprimitive and  $\star=\wr$. In the rest of the section, we focus on the  $\star$-embedding problem for primitive~$G$ and $\star=\uparrow$.

A \emph{cartesian decomposition} of~$\Omega$ is defined in ~\cite{Praeger2018} as a finite set  $\cP=\{P_1,\ldots,P_k\}$ of partitions of $\Omega$ such that $|P_i|\ge 2$ for each~$i$ and $|\Delta_1\cap\cdots\cap\Delta_k|=1$ for each $\Delta_1\in P_1 ,\ldots,\Delta_k\in P_k$. A cartesian decomposition $\cP$ is said to be \emph{trivial} if $\cP$ contains only one partition, namely, the partition into singletons, and $\cP$ is said to be \emph{homogeneous} if the number $|P_i|$ does not depend on~$i=1,\ldots,k$. Every partition~$\pi$ of $\cP$ defines a cartesian decomposition $\cP_\pi$ consisting of  the meets $P_{i_1}\wedge\cdots\wedge\, P_{i_\ell}$, where $\{P_{i_1},\ldots, P_{i_\ell}\}$ is a class of~$\pi$.

A group~$G\le\sym(\Omega)$ \emph{preserves} (respectively, \emph{stabilizes}) the cartesian decomposition $\cP$ if any element of $G$ permutes the $P_i$ (respectively, leaves each $P_i$ fixed). In this case, we say that
$\cP$  is \emph{maximal} for $G$ if $\cP=\cQ_\pi$ for  no cartesian decomposition~$\cQ$ preserved (respectively, stabilized) by~$G$ and no nontrivial partition $\pi$ of~$\cQ$. Note that if $G$ preserves~$\cP$ and the action of $G$ on~$\cP$ is transitive, then~$\cP$ is homogeneous.  Furthermore, if $G$ stabilizes a nontrivial ~$\cP$, then $G$ cannot be primitive.

A natural example of cartesian decomposition comes from the wreath product $G=K\uparrow L$, where as before $K\le\sym(\Delta)$ and $L\le\sym(\Gamma)$. The underlying set of~$G$ is equal to $\Delta^k$, where $k=|\Gamma|$, and one can define a partition $P_i$ ($i=1,\ldots,k$) with $|\Delta|$ classes of the form
$$
\{(\delta_1,\ldots,\delta_k)\in\Delta^k:\ \delta_i\text{ is a fixed element of } \Delta\}.
$$
The partitions $P_1,\ldots, P_k$ form a cartesian decomposition~$\cP$ of~$\Omega$, which is preserved by~$G$ and stabilized by~$K^k$; we say that $\cP$ is a~\emph{standard} cartesian decomposition for~$G$. Clearly, it can be found efficiently for any given~$K$ and~$L$. Well-known properties of a wreath product with the product action~\cite{KaluK1979} imply that if $G$ is primitive, then the standard cartesian decomposition  (a) is homogeneous, and (b)~is maximal (among those that are preserved by $G$) if and only if $K$ is basic.

\lmml{280123s}
Let $G\le\sym(\Omega)$ be a primitive group. Then $G$ is non-basic if and only if $G$ preserves a nontrivial homogeneous cartesian decomposition of~$\Omega$. Moreover, given such a  decomposition, an embedding of $G$ to a wreath product with product action can be found efficiently.
\elmm
\prf
Let $G$ be non-basic. Then there is an embedding of~$G$ to a group~$K\uparrow L$ for some~$K\le\sym(\Delta)$ and~$L\le\sym(\Gamma)$, such that $|\Delta|<|\Omega|$ and $|\Gamma|<|\Omega|$. Denote by $f$ the corresponding bijection from $\Omega$ to $\Delta^{|\Gamma|}$. Then $G$ preserves a homogeneous nontrivial cartesian decomposition $f^{-1}(\cP)$, where $\cP$ is the standard cartesian decomposition for $K\uparrow L$.

Let $G$ preserve a nontrivial homogeneous cartesian decomposition~$P_1,\ldots,P_k$ of~$\Omega$. Denote by $L$ (respectively, $K$) the permutation group induced by the action of~$G$ (respectively, the stabilizer of~$P_1$  in $G$) on the set $\Gamma=\{P_1,\ldots,P_k\}$ (respectively, $\Delta=P_1$). Following the proof of \cite[Theorem~5.13]{Praeger2018}, one can efficiently identify each $P_i$ with~$\Delta$. Then the bijection~$f$ from $\Delta^k= P_1\times\ldots\times P_k$ onto~$\Omega$ taking the cartesian product $\Delta_1\times\ldots\times\Delta_k\in P_1\times\ldots\times P_k$ to the unique point in $\Delta_1\cap\ldots\cap \Delta_k$ can be found efficiently. Now the bijection~$f^{-1}$ moves $G$ to a subgroup of $K\uparrow L$.
\eprf

\thrml{090123b}
Let $G$ be a permutation group of degree~$n$, and $\star\in\{\times,\wr,\uparrow\}$. Assume that $G$ is imprimitive if $\star=\wr$, and primitive if $\star=\uparrow$. Then the  $\star$-embedding problem for~$G$ can be solved in time $\poly(n)$.
\ethrm
\prf
Using standard permutation group algorithms \cite[Section~3.1]{Sere2002}, one can solve the  $\star$-embedding problem for~$G\le\sym(\Omega)$ if $\star=\times$ or $\wr$. Assume that $\star=\uparrow$ and $G$ is primitive. Then (again  by standard permutation group algorithms), one can find in time $\poly(n)$ the socle $S=\soc(G)$ of $G$ and test whether or not $S$ is abelian.

In the abelian case, the required statement can be proved in almost the same way as was done in~\cite[Section~5.1]{EvdP2001b} for solvable groups. Indeed,  in this case, the group~$S$ is elementary abelian of order $n=p^k$ and can naturally be identified with~$\Omega$, which therefore can be treated as a linear space over the field of order~$p$. The procedure \textbf{BLOCK} described in the cited paper, efficiently finds a minimal subspace $\Delta\subseteq \Omega$ so that $\Omega$ is the direct sum of the subspaces belonging to the set $\Gamma=\{\Delta^g:\ g\in G\}$.  Now the required embedding of~$G$ exists only if $\Delta\ne\Omega$, and then as $L$ and $K$ one can take the group $G^\Gamma$ and the restriction of its   stabilizer of $\Delta$ (as a point) to~$\Delta$ (as a set).

Let $S$ be nonabelian. Then $S$ is a direct product of  pairwise isomorphic nonabelian simple groups. We need two auxiliary statements.\smallskip

{\bf Claim 1. }{\it There is at most one  maximal nontrivial cartesian decomposition  $\cP$ stabilized by~$S$. Moreover, one can test  in time $\poly(n)$ whether $\cP$ does exist, and if so, then find it within the same time. }

\prf
We will show that up to the language (in fact, the language of coherent configurations, see \cite{CP2019}) this claim is an almost direct consequence of results in~\cite{Chen2021,KN2009}. We start by noting that the cartesian decompositions stabilized by $S$ are exactly the tensor decompositions of the coherent configuration $\cX$ associated with~$S$ (see \cite{Chen2021} for details).  Thus, in view of \cite[Theorem~1]{Chen2021}, the cartesian decompositions stabilized by $S$ are in a 1-1 (efficiently computable) correspondence with the  cartesian decompositions of the coherent configuration $\cX$ itself.  Moreover, if every subdegree of~$S$ is at least~$2$, i.e.,~$\cX$ is thick in terms of~\cite{Chen2021},  then there is at most one  maximal nontrivial cartesian decomposition~$\cP$ of $\cX$ \cite[Theorem~2]{Chen2021}. The polynomial-time algorithm  in  \cite[Lemma~13]{Chen2021} enables us to find a certificate that $\cX$ has only the trivial cartesian decomposition, or to construct~$\cP$.

Assume that at least one (nontrivial) subdegree of $S$ is equal to~$1$. Since the union of singleton orbits of a one point stabilizer of~$S$ is a block of the primitive group $G$,  this union is the whole set $\Omega$ and the group $S$ is regular. In this case, the coherent configuration~$\cX$ is also regular, $S=\aut(\cX)$, and from the above mentioned~\cite[Theorem~1]{Chen2021}, it follows that the  cartesian decompositions of~$\cX$ are in a 1-1 correspondence with the direct decompositions of the group~$S$ itself. If this group is simple, then $S$ stabilizes only the trivial cartesian decomposition. Otherwise,  the decomposition of~$S$ into the direct product of pairwise isomorphic (nonabelian) simple groups gives the  maximal nontrivial cartesian decomposition  $\cP$ stabilized by~$S$. It remains to note that $\cP$ can be found efficiently by the main  algorithm in~\cite{KN2009}.
\eprf

{\bf Claim 2.}{\it\  Assume that $G$ is non-basic. Then $G$ preserves a nontrivial homogeneous cartesian  decomposition of the form~$\cP_\pi$ for some partition~$\pi$ of the cartesian decomposition~$\cP$ from Claim~1.}

\prf
Since $G$ is non-basic, we may assume that $G\le K\uparrow L$, where $K$ is basic primitive and $L$ is transitive. Denote by~$Q$ the corresponding standard cartesian decomposition (Lemma~\ref{280123s}).  We may also assume that $K$ is the permutation group induced by the action on $P\in Q$ of the stabilizer of $P$ in $G$. Then in virtue of~\cite[Theorem~4.7]{Cam1999} (and the remark after it), the socle~$S$ of~$G$ is a subgroup of the base group of the wreath product~ $K\uparrow L$. It follows that~$S$ stabilizes $Q$. Thus, by Claim~1, there exists the unique maximal nontrivial cartesian decomposition  $\cP$ stabilized by~$S$ and $Q=\cP_\pi$ for some partition~$\pi$ of~$\cP$. Since $G$ acts transitively on~$\cP$, the decomposition $Q$ is homogeneous.
\eprf

Let us complete the proof. By Lemma~\ref{280123s}, it suffices to test whether~$G$ preserves a nontrivial cartesian decomposition and, if so, find it efficiently. Applying the algorithm of Claim~1, we test efficiently  whether ~$S$ stabilizes a nontrivial cartesian decomposition. If not, then $G$  cannot preserve  a nontrivial cartesian decomposition (see Claim~2). Otherwise, we efficiently find the cartesian decomposition $\cP$ from Claim~1. By  Claim~2,  all we need is  to test whether there exist a partition~$\pi$ of~$\cP$, such that~$\cP_\pi$ is a nontrivial homogeneous cartesian decomposition preserved by~$G$ and, if so, find it efficiently. Since $\cP$ is nontrivial, we have $|\cP|\le \log n$ and the power set~$2^\cP$ has cardinality at most $2^{\log n}=n$. Furthermore, the cartesian decompositions $\cP_\pi$  preserved by~$G$ are in one-to-one correspondence with those  subsets $Q\subseteq \cP$ for which
$$
\{Q^g:\ g\in G\}\text{ is a homogeneous partition of } \cP.
$$
 Since this condition can be tested only for the generators~$g$ of $G$, we are done.
\eprf

\section{Proof of Theorem~\ref{040123b}}\label{040123b0}

We deduce Theorem~\ref{040123b}  from a more general statement valid for any complete class of groups. A class of (abstract) groups is said to be \emph{complete} if it is closed with respect to taking subgroups, quotients, and extensions \cite[Definition~11.3]{Wie1964}. Any complete class is obviously closed with respect to direct and wreath product and taking sections. The class of all permutation groups of degree at most~$n$ that belong to~$\fK$ is denoted by~$\fK_n$.

\thrml{040123az}
Let $m,n\in\mN$, $m\ge 3$, and $\fK$ a complete class of groups. Then
\nmrt
\tm{i} $\fK_n$ is closed with respect to taking the $m$-closure if and only if $\fK_n$ contains the $m$-closure of every primitive basic group in~$\fK_n$,
\tm{ii} the $m$-closure of any  group in~$\fK_n$ can be found in time $\poly(n)$ by 
accessing oracles for finding the $m$-closure of every  primitive basic group in $\fK_n$ and the relative $m$-closure of every group in~$\fK_n$ with respect to any group in~$\fK_n$.
\enmrt
\ethrm
\prf
The ``only if'' part of statement~(i) is obvious. To prove the ``if''  part and statement (ii), we present a more or less standard recursive algorithm for finding the $m$-closure~$G^{(m)}$ of a  group $G\in\fK_n$. At each step we will verify that  $G^{(m)}\in\fK_n$.

Depending on whether $G\le\sym(\Omega)$ is intransitive, imprimitive, or primitive, we set $\star=\times$, $\wr$, or~$\uparrow$, respectively.  Solving the  $\star$-embedding problem for $G$ by Theorem~\ref{090123b}, one can can test  in time $\poly(n)$  whether there exists an embedding of $G$ to $K\star L$ for some sections
$$
K\le\sym(\Delta)\qaq L\le\sym(\Gamma)
$$
of $G$, such that the numbers $n_K=|\Delta|$ and $n_L=|\Gamma|$ are less than $n=|\Omega|$, and if so, then find the embedding explicitly. If there is no such embedding, then $G$ is primitive basic,  belongs to $\fK_n$, and   the $m$-closure of~$G$  can be found for the cost of one call of the corresponding oracle.

Assume that  $G$ is not  primitive basic and we are given a bijection~$f$ from $\Omega$ to the underlying set of~$K\star L$, such that equality~\eqref{010223e} holds. Since  $f^{-1}G^{(m)}f=(f^{-1}Gf)^{(m)}$, we may also assume that
$$
G\le K\star L.
$$

Note that $K\in\fK_{n_K}$ and $L\in\fK_{n_L}$, because the class $\fK$ is complete. Applying the algorithm recursively to $K$ and $L$, we find the groups $K^{(m)}$ and $L^{(m)}$ in time $\poly(n_K)$ and $\poly(n_L)$, respectively, and then the group $K^{(m)}\star L^{(m)}$ in time $\poly(n)$. By induction, $K^{(m)}\in\fK_{n_K}$ and $L^{(m)}\in\fK_{n_L}$, whence $K^{(m)}\star L^{(m)}\in\fK_n$. On the other hand, by Theorem~\ref{090123a}, we have
$$
G^{(m)}\le (K\star L)^{(m)}\le  K^{(m)}\star L^{(m)}.
$$
Thus, $G^{(m)}\in\fK_n$. Accessing  (one time) the oracle for finding the relative $m$-closure of $G$ with respect to $K^{(m)}\star L^{(m)}$, we finally get the group $G^{(m)}$.

It remains to estimate the number of the oracles calls. Each recursive call divides the problem for a group of degree~$n$ to the same problem for a group $K$ of degree~$n_K$ and for a group $L$ of degree~$n_L$. Moreover,
$$
n=\css
n_K+n_L &\text{if $\star=\times$,}\\
n_K\cdot n_L &\text{if $\star=\wr$,}\\
{n_K}^{n_L} &\text{if $\star=\uparrow$.}\\
\ecss
$$
Thus the total number of recursive calls and hence the number of accessing oracles is at most $n$.
\eprf

An obstacle in proving Theorem~\ref{040123az} for $m=2$ lies in the exceptional case of Theorem~\ref{090123a}. Indeed, assume that the class $\fK$ does not contain all groups. Then it cannot contain symmetric groups of arbitrarily large degree. However the $2$-closure of any two-transitive group of degree $n$ coincides with~$\sym(n)$. Therefore $\fK$ cannot also contain two-transitive groups of sufficiently large degree. It seems that this restricts the class $\fK$ essentially.

\rmrkl{040123a}
In fact, the proof of Theorem~{\rm\ref{040123az}} shows that the following  weakened version of this theorem holds true: both statements of Theorem~{\rm\ref{040123az}} remain valid for $m=2$ if ``primitive basic groups'' are replaced with ``primitive groups''.
\ermrk

{\bf Proof of Theorem~\ref{040123b}.} Denote by $\fK$ the class of all solvable groups. This class is obviously complete. Moreover, the relative $m$-closure of any group of~$\fK_n$ with respect to any other group from~$\fK_n$ can be found in time~$\poly(n)$ in view of~\cite[Corollary~3.6]{BabaL1983} (see also \cite[Section 6.2]{Ponomarenko2020b}). By Theorem~\ref{040123az}, it suffices to verify that the $3$-closure of a primitive basic group $G\in\fK_n$ can be found in time $\poly(n)$; indeed, if $m>3$, then $G^{(m)}\le G^{(3)}$ can be found as the relative $m$-closure of $G$ with respect to~$G^{(3)}$.

First, suppose that a point stabilizer $H$ of~$G$ has a regular orbit. Then $G$ is $3$-closed  by \cite[Corollary~2.5]{OBrien2020}, and there is nothing to do, because $G=G^{(3)}$. Now, if  the group $H$ has no regular orbits and $n$ is  sufficiently large, then the number $n=q$ is a prime power and $H\le\GaL(1,q)$, see~\cite[Corollary~3.3]{Yang2020}. In this case, $H=H^{(2)}$ by \cite[Proposition~3.1.1]{Xu2011d}  and again $G=G^{(3)}$.  In the remaining case, the degree of $G$ is bounded by an absolute constant, say $N$, and the group~$G^{(3)}$ can be found by inspecting all permutations of~$\sym(N)$.

\end{document}